\documentclass[11pt]{article}
\usepackage{amsmath,amssymb} 

\begin{document}

\title{Seifert fibered surgery on Montesinos knots}

\author{Ying-Qing Wu} 
\date{}
\maketitle



\begin{abstract}
  Exceptional Dehn surgeries on arborescent knots have been classified
  except for Seifert fibered surgeries on Montesinos knots of length
  3.  There are infinitely many of them as it is known that $4n+6$ and
  $4n+7$ surgeries on a $(-2,3,2n+1)$ pretzel knot are Seifert
  fibered.  It will be shown that there are only finitely many others.
  A list of 20 surgeries will be given and proved to be Seifert
  fibered.  We conjecture that this is a complete list.
\end{abstract}

\newcommand{\proof}{\noindent {\bf Proof.} }
\newcommand{\qed}{\quad $\Box$}
\newtheorem{thm}{Theorem}[section]
\newtheorem{prop}[thm]{Proposition} 
\newtheorem{lemma}[thm]{Lemma} 
\newtheorem{sublemma}[thm]{Sublemma} 
\newtheorem{cor}[thm]{Corollary} 
\newtheorem{defn}[thm]{Definition} 
\newtheorem{convention}[thm]{Convention} 
\newtheorem{notation}[thm]{Notation} 
\newtheorem{qtn}[thm]{Question} 
\newtheorem{example}[thm]{Example} 
\newtheorem{remark}[thm]{Remark} 
\newtheorem{conj}[thm]{Conjecture} 
\newtheorem{prob}[thm]{Problem} 
\newtheorem{rem}[thm]{Remark} 

\newcommand{\bdd}{\partial}
\newcommand{\Int}{{\rm Int}}
\newcommand{\wind}{{\rm wind}}
\newcommand{\wrap}{{\rm wrap}}
\renewcommand{\a}{\alpha}
\renewcommand{\b}{\beta}

\input epsf.tex

\section{Introduction}

A Dehn surgery on a hyperbolic knot is {\it exceptional\/} if it
is reducible, toroidal, or Seifert fibered.  By Perelman's work, all
other surgeries are hyperbolic.  For knots in $S^3$, by exceptional
surgery we shall always mean {\it nontrivial\/} exceptional surgery.

Given an arborescent knot, we would like to know exactly which
surgeries are exceptional.  We divide arborescent knots into three
types.  An arborescent knot is of type I if it has no Conway sphere,
so it is either a 2-bridge knot or a Montesinos knot of length 3.
A type II knot has a Conway sphere cutting it into two tangles, each
of which is the sum of two nontrivial rational tangles, with one of them
of slope $1/2$.  All others are of type III.  In [Wu1] it was shown
that all nontrivial surgeries on type III arborescent knots are Haken
and hyperbolic, and all nontrivial surgeries on type II knots are
laminar.  In [Wu2] it was further shown that there are exactly three
type II knots admitting exceptional surgery, each of which admits
exactly one exceptional surgery, producing a toroidal manifold.  For
type I knots, Brittenham and the author determined exceptional
surgeries on 2-bridge knots [BW], toroidal surgeries on Montesinos
knots have been classified in [Wu3], and it is known that there is no
reducible surgery on hyperbolic arborescent knots [Wu1].

It remains to determine small Seifert fibered surgeries on hyperbolic
Montesinos knots of length 3, which is also the set of all Seifert
fibered surgeries because by Ichihara and Jong [IJ1] the only {\it
  toroidal\/} Seifert fibered surgery on Montesinos knots is the 0
surgery on the trefoil knot and hence there is no large Seifert
fibered surgery on hyperbolic Montesinos knots.  For the special case
of finite surgeries on Montesinos knots, the classification has been
done by Ichihara and Jong [IJ2].  See also [FIKMS].

In general, let $K = K(p_1/q_1,\, p_2/q_2,\, p_3/q_3)$ be a hyperbolic
Montesinos knot of length 3, and assume that it admits a nontrivial
Seifert fibered surgery $K(r)$.  Using immersed surfaces, it was shown
in [Wu4] that we must have $\frac 1{q_1-1} + \frac 1{q_2-1} + \frac
1{q_3-1} \leq 1$, hence up to relabeling we have $|q_1|=2$, or
$|q_1|=|q_2|=3$, or $(|q_1|, |q_2|, |q_3|) = (3,4,5)$.  In [Wu5] we
studies persistently laminar branched surfaces in knot complements,
and obtained restrictions on the $|p_i|$.  More explicitly, if $K$
above is a pretzel knot $K(1/q_1,\, 1/q_2,\, 1/q_3,\, n)$ then either
(i) $n=0$, or (ii) $n=-1$ and $q_i>0$, and if $K$ is not a pretzel knot
then it is either (iii) $K(2/3,\, 1/3,\, 2/5)$, or (iv) $K(1/2,\,
1/3,\, 2/(2a + 1))$ with $a \in \{3, 4, 5, 6\}$, or (v) $K(1/2,\,
2/5,\, 1/q)$ for some odd $q\geq 3$.

There are still infinitely many knots among the above, for example it
includes all $(2, q_2, q_3)$ pretzel knots.  In [Wu6] we studies
exceptional surgeries on tubed Montesinos knots.  These are the knots
in solid tori obtained by tubing Montesinos tangles in some specific
ways.  By embedding the solid tori into $S^3$, we see that these knots
are closely related to Montesinos knots in $S^3$.  With this method it
was shown that there are indeed infinitely many Seifert fibered
surgeries on Montesinos knot of length 3, that is, each $(-2,3,2n+1)$
pretzel knot in $S^3$ admits at least two Seifert fibered surgeries,
of slopes $4n+6$ and $4n+7$, respectively.  On the other hand, using
the classification theorem in [Wu6], we will show that there are only
finitely many other Seifert fibered surgeries on these knots.  See
Theorem 2.3 below.

A few other surgeries on Montesinos knots are known to be Seifert
fibered.  There is that well known $(-2,3,7)$ pretzel knot, on which
17, 18 and 19 surgeries are Seifert fibered [FS].  Hyun-Jong Song
showed that surgery on $(-3,3,3)$ with slope $1$ is Seifert fibered,
and Mattman, Miyazaki and Motegi [MMM] showed that surgery on
$(-3,3,5)$ of slope $1$ is also Seifert fibered.  More examples will
be given in Table 3.1 below.  We conjecture that this is a complete
list.

\section{A finiteness theorem}

Consider the knots $K_n = K(-1/2,\, 1/3,\, 1/2n+1)$.  By [Wu6,
Corollary 2.3], $K_n(r_n)$ is small Seifert fibered for $r_n = 6+4n$
and $7+4n$, except that $K_2(15)$ is reducible.  If $K'$ is the mirror
image of $K$, then an orientation reversing homeomorphism of $S^3$
induces an orientation reversing homeomorphism from $K(r)$ to $K'(-r)$.
We consider $(K,r)$ and $(K',-r)$ as equivalent.  Theorem 2.3 below
shows that there are only finitely many other small Seifert fibered
surgeries on length 3 Montesinos knots.

More generally, consider a two component link $L = K' \cup K''$ with
$K''$ a trivial component.  Let $V$ be the solid torus $S^3 - \Int
N(K'')$ , and let $(V, K', r)$ be the manifold obtained by $r$ surgery
on $K'$ in $V$.  Denote by $K_m$ the knot obtained from $K'$ by $m$
right-hand full twists on $K''$.  

\begin{lemma} Suppose $L = K' \cup K''$ is a two component hyperbolic
  link in $S^3$ with $K''$ a trivial loop.  Then there is a finite
  collection $C$ of $(m, r_n)$, where $r_n$ is a slope on $\bdd
  N(K')$, such that if $K_m(r_n)$ is non-hyperbolic, then either (i)
  $K_m$ is non-hyperbolic, or (ii) $(V, K', r_n)$ is nonhyperbolic, or
  (iii) $(m, r_n) \in C$.
\end{lemma}

\proof This is well known and follows immediately from the
$2\pi$-theorem of Gromov and Thurston.  By the $2\pi$-theorem there is
a finite set $C_i$ of slopes on each cusp $T_i$ of $S^3 - \Int N(L)$,
such that if $r_i \not \in C_i$ for $i=1,2$ then $L(r_1, r_2)$ is
hyperbolic.  Let $\hat C$ be the collection of all slopes $(r_1, r_2)$
such that $L(r_1, r_2)$ is non-hyperbolic, and let $C'_i$ be the set
of slopes $r$ on $T_i$ such that $r$ filling on $T_i$ is
non-hyperbolic.  If for some $r_1$ there are infinitely many $r_2$
such that $(r_1, r_2) \in \hat C$ then $r_1 \in C'_1$.  Similarly if
there are infinitely many $r_1$ with $(r_1, r_2) \in C$ then $r_2 \in
\hat C'_2$.  Thus if we denote by $\hat C_i = \{(r_1, r_2) \; | \; r_i
\in C'_i\}$ then $C = \hat C - \hat C_1 \cup \hat C_2$ is finite.
Restricting the above to the set of slopes with $r_1$ of type $1/m$
gives the required result.  \qed

\medskip 

The following result was conjectured by Gordon and proved by Lackenby
and Meyerhoff [LM].  It will be referred to as the $8$-Theorem below.

\bigskip

\noindent
{\bf The 8-Theorem (Lackenby-Meyerhoff)}  \quad {\em
If $M$ is a hyperbolic manifold and $r_1, r_2$ are two exceptional 
slopes on a torus component of $\bdd M$, then $\Delta(r_1, r_2) \leq 8$.
}

\bigskip

Consider the links $L = K' \cup K''$ in Figure 2.1 below, where $K''$
is the trivial knot, and $K'$ is a in the solid torus $V = S^3 - \Int
N(K'')$ obtained by adding two strings to a Montesinos tangle
$T(p_1/q_1,\, p_2/q_2)$, as shown in Figure 2.1(a)-(b).  These are
called {\it tubed Montesinos knots\/} in [Wu6].  Denote by
$K^0(p_1/q_1,\ p_2/q_2)$ the knot in Figure 2.1(a), and by
$K^1(p_1/q_1,\ p_2/q_2)$ the one in Figure 2.1(b).  We always assume
that $K'$ is a knot in $V$, so $q_1, q_2$ are not both even.  Denote
by $(V, K', r)$ the manifold obtained by $r$ surgery on $K'$ in $V$.
A knot in $V$ is considered to be equivalent to its mirror image in
the lemma below.

\bigskip
\leavevmode

\centerline{\epsfbox{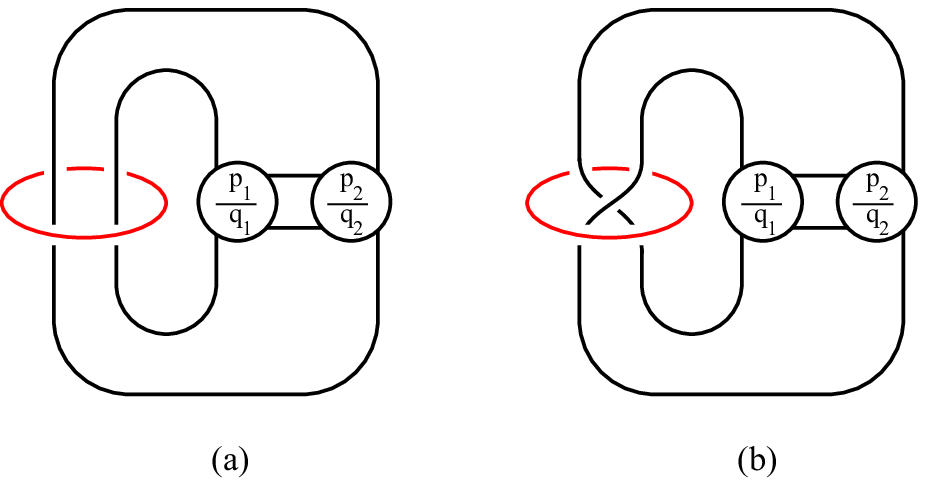}}
\bigskip
\centerline{Figure 2.1}
\bigskip

\begin{lemma} \label{lemma:2.2} Suppose $q_i \geq 2$, $r$ is a
  nontrivial slope, and $(V, K', r)$ is non-hyperbolic.  Then $(K',
  r)$ is equivalent to one of the following pairs.  The surgery is
  small Seifert fibered for $r=7$ in (2), and toroidal otherwise.

  (1) $K = K^a(1/q_1,\, 1/q_2)$, $|q_i|\geq 2$, $a=0,1$, and $r$ is
  the pretzel slope.

  (2) $K = K^1(-1/2,\, 1/3)$, $r = 6,7,8$.
\end{lemma}

\proof Exceptional surgeries for all tubed Montesinos knots in solid
torus have been classified in [Wu6, Theorem 5.5].  This lemma is the
above theorem applied to the case that $K'$ is a tubed Montesinos knot
of length 2 in $V$.  \qed

\begin{thm} Besides the $4n+6$ and $4n+7$ surgeries on the
  $(-2,3,2n+1)$ pretzel knots, there are only finitely many small
  Seifert fibered surgeries on hyperbolic Montesinos knots $K$ of
  length 3.  Moreover, $K$ is equivalent to one of the following.

  (1) $(q_1, q_2, q_3)$ pretzel knot, $|q_1|\leq |q_2| \leq |q_3| \leq
  17$, and either $|q_1| = 2$ or $|q_1| = |q_2| = 3$.

  (2) $(3,\,3,\,2n,\,-1)$ pretzel knot, $2 \leq n \leq 8$.

  (3) $K(-1/2,\, 2/5,\, 1/(2n+1))$ for some $n>0$.

  (4) Ten individual knots: $(3,\, \pm 4,\, \pm 5)$
  pretzel knots, $(3,\, 4,\, 5,\, -1)$ pretzel knot, $K(-2/3,\, 1/3,\,
  2/5)$, and $K(-1/2,\, 1/3,\, 2/(2a+1))$ for $a=3,4,5,6$. 
\end{thm}

\proof By [Wu5], if $K$ is a hyperbolic Montesinos knot of length 3
and $K(r)$ is atoroidal and Seifert fibered, then $K$ is one of the
following knots.

(a) $K = (q_1, q_2, q_3)$ pretzel knot, and either $|q_1| = 2$,
or $|q_1| = |q_2| = 3$, or $(|q_1|, |q_2|, |q_3|) = (3,4,5)$;

(b) $K = (q_1, q_2, q_3, -1)$ pretzel knot with $q_i \geq 3$, and either
$q_1 = q_2 = 3$ or $(q_1, q_2, q_3) = (3,4,5)$;

(c) $K = K(-2/3,\, 1/3,\, 2/5)$;

(d) $K = K(1/2,\, 1/3,\, 2/(2a + 1))$ and $a \in \{3, 4, 5, 6\}$;

(e) $K = K(-1/2,\, 2/5,\, 1/q)$ for some $q \geq 3$ odd.

\medskip

Besides the 10 individual knots listed in (5), there are several
infinite families of knots among the above.  We divide these into four
cases as follows.

\medskip 

Case 1.  {\em $K$ is a $(q_1, q_2, q_3)$ pretzel knot with $2 = |q_1|
  \leq |q_2| \leq |q_3|$.}

\medskip 
 
Up to equivalence we may assume $q_1 = -2$.  Let $L = K_1 \cup K_2$ be
the link in Figure 2.1(b), where $p_1/q_1 = -1/2$, $p_2/q_2 = 1/q_2$,
and $K_1$ denotes the trivial circle in the figure.  Denote by $L(r_1,
r_2)$ the manifold obtained from $L$ by $r_i$ surgery on $K_i$.  By
Kirby Calculus, we see that $L(1/n, r_2) = K(r)$, where $K = K(-1/2,\,
1/q_2,\, 1/(1-2n))$, and $r = r_2 - 4n$.  On the other hand, if we
denote by $V$ the exterior of $K_1$, and put $M = (V,K_2,r_2)$, then
$L(r_1, r_2) = M(r_1)$, the manifold obtained by Dehn filling along
slope $r_1$ on $\bdd M$.

By Lemma 2.2, if $M$ is non-hyperbolic then $q_2 = 3$ and $r_2 =
6,7,8$, in which case $K(r)$ is a $6-4n$, $7-4n$ or $8-4n$ surgery on
$K(-1/2, 1/3, 1/(1-2n))$.  The first two are Seifert fibered and have
been excluded in the statement, while the last one is the pretzel
slope, in which case $K(r)$ is toroidal [Wu3] and hence cannot be
Seifert fibered [IJ1].  Thus we may assume that $M = (V,K_2,r_2)$
is hyperbolic.  It is easy to see that $M(1/0)$ is the $r_2$ surgery
on the $(2,5)$ torus knot and hence is nonhyperbolic.  Therefore by
the 8-Theorem we see that if $M(1/n)$ is non-hyperbolic then $|n|\leq
8$, hence $|q_3| = |1-2n| \leq 17$.

\medskip 

Case 2.  {\em $K$ is a $(q_1, q_2, q_3)$ pretzel knot with $3 = |q_1|
  = |q_2| \leq |q_3|$.}

\medskip 

Similar to the above, we have $L = K_1 \cup K_2$ as in Figure 2.1(a)
or 2.1(b), according to whether $q_3 = 2n$ or $2n+1$.  In this case
$|q_1| = |q_2| = 3$, so $K_2 = K^0(\pm 3, 3)$ or $K^1(\pm 3, 3)$.  By
Lemma 2.2, the only exceptional surgery on $K_2$ in $V$ is the surgery
along the pretzel slope, corresponding to the toroidal surgery along
the pretzel slope of $K$ [Wu3].  For all other nontrivial slopes $r$,
$M = (V,K_2,r)$ is hyperbolic.  Note that when $q_3 = 2n$, $M(1/0)$ is
the $r$-surgery on the connected sum of two trefoil knot and hence is
nonhyperbolic, therefore as above, we see that if $|n|>8$ then $r$
surgery on $K$ is hyperbolic.  When $q_3 = 2n+1$ is old and $q_1 = q_2
= 3$, the knot with $n = -1$ is the trefoil knot, so the above
argument applies and we conclude that $|n-(-1)| \leq 8$ if $K$ admits
a nontrivial small Seifert fibered surgery.  Thus $|q_3| = |2n +1|
\leq 17$.

Now consider the case that $K$ is the $(-3,3,2n+1)$ pretzel knot.  In
this case we need to use a recent result of Boyer, Gordon and Zhang.
Since $K$ has a Seifert surface of genus 1, by [BGZ, Theorem 1.5]
$K(p/q)$ is hyperbolic unless $|p|\leq 3$.  By [Wu5, Theorem 6.6] we
know that the knot complement has a persistently laminar branched
surface with two meridional cusps, hence if $q\neq 1$ then the
lamination is genuine in $K(p/q)$, so by [Br] $K(p/q)$ cannot be
Seifert fibered.  $K$ is the twist knot $6_1$ when $n = 0$, and its
mirror image when $n = -1$, one of which has small Seifert surgery
slopes $1,2,3$ and the other $-1,-2,-3$.  It now follows by the same
argument as above that for $p/q = \pm 1,\pm2, \pm 3$, the surgery
$K(p/q)$ is hyperbolic unless $|q_3| \leq 17$.

\medskip 

Case 3.  {\em $K$ is a $(3,3,-1,q_3)$ pretzel knot with $q_3 \geq
  3$. }

\medskip 

Since $K$ is a knot, $q_3$ must be even, say $q_3 = 2n$, so we have
$n\geq 2$.  As above, let $L = K_1 \cup K_2$, where $K_1$ is trivial
and $K_2$ is a tubed knot $K^0(1/3, -2/3)$ in the solid torus $V = S^3
- \Int N(K_1)$.  Now $r$ surgery on $K$ is equivalent to $r-4n$
surgery on $K_2$ followed by $1/n$ Dehn filling on $\bdd V$, with
respect to the standard meridian-longitude coordinate of $K_1$.  By
Lemma 2.2 there is no exceptional surgery on $K_2$ in $V$, hence
$(V,K_2,r-4n)$ is hyperbolic. Also note that when $n=0$ the knot $K$ is
the connected sum of two trefoils, hence all surgeries are
non-hyperbolic.  It follows that $K(r)$ is hyperbolic unless $n \leq
8$.
 
\medskip 

Case 4.  {\em $K = K(-1/2,\, 2/5,\, 1/(2n+1)$ for some $n \geq 1$.}

\medskip 

Let $L = K_1 \cup K_2$, where $K_1$ is trivial and $K_2 = K^1(-1/2,
2/5)$ in $V = S^3 - \Int N(K_1)$.  $K(r)$ is the same as $r-4n$
surgery on $K_2$ followed by $1/n$ filling on $\bdd V$.  By Lemma 2.2
we see that $(V, K_2, r)$ is always hyperbolic for any nontrivial $r$,
hence by Lemma 2.1  there are only finitely many exceptional surgeries
on the set of hyperbolic knots $K$ as above.  \qed

\medskip

We note that the argument above does not provide a bound for $n$ for
the knots of type $K(-1/2, 2/5, 1/(2n+1))$, although by the theorem
such bound does exist.  However, using computer assistant proof it
seems likely that $n \leq 9$.  See the discussion about {\it
  Snappex\/} after Conjecture 4.1.

\section{Seifert fibered surgeries}

Gordon conjectured that a Seifert fibered surgery on a hyperbolic knot
is an integral surgery.  The following lemma shown that this is
true for most Montesinos knots of length 3.

\begin{lemma} Suppose $K$ is a hyperbolic Montesinos knot of length 3
  such that $K(r)$ is small Seifert fibered and $r$ is a nontrivial
  non-integral slope. Then $K$ is equivalent to either (i) a $(-2,p_2,
  p_3)$ pretzel knot with $3 \leq p_2 \leq p_3 \leq 17$, or (ii) a
  $(3,3,-1,2n)$ knot with $2\leq n \leq 8$, or (iii) the $(3,4,5,-1)$
  pretzel knot.
\end{lemma}

\proof By the proof of [Wu6, Theorem 6.6] $K$ has a persistently
laminar branched surface with two meridional cusps unless it is a
pretzel knot of type $(p_1, p_2, p_3, -1)$ with $p_i>1$.  Such
branched surface becomes genuinely laminar after nonintegral surgery
because the component containing the Dehn filling solid torus is a
solid torus with cusps intersecting a meridian disk at least 4 times,
hence by [Br] the surgered manifold cannot be a small Seifert fibered
manifold.  The result follows by comparing this with the list of knots
in Theorem 2.3.  \qed

\medskip

The lemma can be used in searching for Seifert fibered surgeries.  We
may now use Snappy [CDW] to check surgeries on the list of knots in
Theorem 2.3.  $K(r)$ is likely to be Seifert fibered if the program
gives nearly zero volume.  Since most of those knots in the list are
strongly invertible, one can then try to use the Montesinos trick to
show that the manifold is indeed Seifert fibered.  The following table
gives the list of Seifert fibered surgeries on hyperbolic Montesinos
knots of length 3, where $M(r_1, r_2, r_3)$ denotes the closed
3-manifold which is the double branched cover of $S^3$ with branch set
a Montesinos link $K(r_1, r_2, r_3)$.  It is well known that all such
$M(r_1, r_2, r_3)$ are small Seifert fibered.

$$
\begin{array}{llll}
    & K  & r & K(r) \\
\hline 
(1) &  K(-1/2,\, 1/3,\, 1/2n+1) & r = 4n+6 & M(1/2,\ -1/4,\ 2/2n-5) \\
&  &r = 4n+7 & M(-1/3,\, 3/5,\, 1/(n-2)) \\
(2) & K(-1/2,\, 1/3,\, 1/7) &r = 17 & M(-1/2,\, 1/3,\, -2/5) \\
(3) & K(-1/2,\, 1/3,\, 2/5) &r = 3 &  M(-2/15,\, 1/2,\, -1/3)\\ 
&  &r = 4 &  M(-2/7,\, 1/2,\, -1/6)\\
&  &r = 5 & M(3/5,\, -1/3,\, -1/5)\\
(4) & K(-1/2,\, 1/5,\, 2/5) &r = 7 & M(3/4,\, -2/5,\, -1/4)\\
&  &r = 8 & M(-1/5,\, 1/2,\, -2/9) \\
(5) & K(-1/2,\, 1/7,\, 2/5) &r = 11 & M(-1/4,\, -2/7,\, 2/3)  \\
(6) & K(-1/2,\, 1/3,\, 2/7) &r = -1 & M(-3/4,\, 1/3,\, 3/8)\\
&  &r = 0 & M(1/5,\, 3/10,\,-1/2)\\
&  &r = 1 & M(1/2,\, -2/3,\, 3/19)\\
(7) & K(-1/2,\, 1/3,\, 2/9) & r = 2 & M(-3/8,\, -3/2,\, -1/4) \\
&  &r = 3 & M(8/11, -1/2, -1/5) \\
&  &r = 4 & M(-3/20, -1/2, 2/3) \\
(8) & K(-1/2,\, 1/3,\, 2/11) &r = -2 & M(2/7,\, 2/5,\, -2/3) \\
&  &r =  -1 & M(2/9, 2/7, -1/2) \\
(9) & K(-1/3,\, 1/3,\, 1/4) &r = 1 & M(-1/2,\, 1/5,\, 2/7) \\
(10) & K(-1/3,\, 1/3,\, 1/6) &r = 1 & M(-1/2,\, 1/3,\, 2/13) \\ 
(11)  & K(-1/3,\, 1/3,\, 1/3) &r = 1 & M(1/2, -1/5, -2/7)  \\
(12) & K(-1/3,\, 1/3,\, 1/5) &r = 1 & M(-1/3, -1/4, 3/5)  \\
(13) & K(-2/3,\, 1/3,\, 2/5) &r = -5 & M(2/5,\, 2/5,\, -3/4) \\
\end{array}
$$
\centerline{Table 3.1 \qquad Seifert fibered surgeries}

\bigskip

\begin{thm} For each knot $K$ and slope $r$ in the table, $r$ surgery
  on $K$ produces a Seifert fibered manifold $K(r)$ as shown in the
  table. 
\end{thm}

\proof (1) is given in [Wu6, Theorem 5.5].  (2) is well known, see for
example [CGLS].  Most of the others can be proved using the Montesinos
trick.
 
Consider a strongly invertible knot $K$ in $S^3$ with axis $X$
intersecting $K$ twice.  $\pi$-rotation along $X$ gives a quotient map
$\rho: (S^3, X, K) \to (\bar S, \bar X, \bar K)$, where $\bar S$ is a
3-sphere, $\bar X$ a trivial circle, and $\bar K$ an arc with its two
endpoints on $\bar X$.  For example, when $K$ is the knot $K(-1/2,\
1/3,\ 2/5)$ in Figure 3.1(1), the pair $\bar X$ and $\bar K$ are shown
in Figure 3.1(2).  The quotient of $N(K)$ is a 3-ball $\bar N$ in
$\bar S$, drawn as a thick arc in Figure 3.1(3).  Shrinking $\bar N$
to a round ball gives Figure 3.1(4).  Put $\alpha = \bar X \cap \bar
N$.  We may consider $(\bar N, \alpha)$ as a rational tangle of slope
$\infty$, and set up coordinates so that a longitude on $\bdd N(K)$
projects to a curve of slope $0$ on $\bdd \bar N$, which is considered
as a pillow case with the four points $X \cap \bdd \bar N$ as cone
point, so every essential simple closed curve on $\bdd \bar N$ has a
slope; see [HT].  The Montesinos trick [Mon] says that $K(r)$ can be
obtained by replacing $(\bar N, \alpha)$ with a rational tangle of
slope $-r$ to obtain a link $L[-r]$ in $\bar S$, then taking the
double branched cover of $\bar S$ along $L[-r]$.  More generally, if
$\bar N$ is deformed so that $\alpha$ is the $\infty$ tangle and the
longitude projects to a curve of slope $r_0$ on $\bdd \bar N$ then a
curve of slope $r$ on $\bdd N(K)$ projects to a curve of slope $r_0 -
r$ on $\bdd \bar N$, hence $K(r)$ is the double branched cover of
$\bar S$ along $L[r_0 - r]$.  Thus if $L[r_0-r]$ is a Montesinos link
$K(a,b,c)$ then $K(r)$ is the Seifert fibered manifold $M(a,b,c)$.

Continue with the example above.  We can simplify Figure 3.1(4) to that
of Figure 3.1(5).  To determine the framing, consider the bounded
checkboard surface $F$ for the diagram in Figure 3.1(1).  It is easy to
check that $\bdd F$ has slope $6$ on $\bdd N(K)$.  The quotient of $F$
is a disk $\bar F$ which deforms to a disk in Figure 3.1(3) whose
intersection with $\bdd \bar N$ is an arc of slope $0$.  In other
words, the longitude of $N(K)$ projects to a curve of slope $r_0$ on
$\bdd N(K)$ such that a curve of slope $r = 6$ projects to a curve of
slope $r_0 - r = 0$, hence $r_0 = 6$.  Thus $K(3)$ is the double
branched cover of $L[r_0-3] = L[3]$, shown in Figure 3.1(6).  It can be
deformed to that in Figure 3.1(7) and then further to Figure 3.1(8), which
is a Montesinos knot $K(-2/15,\; 1/2,\; -1/3)$.  Hence $K(3) =
M(-2/15,\; 1/2,\; -1/3)$.  Similarly, $K(4)$ is the double branched
cover of $L[r_0-4] = L[2]$ shown in Figure 3.1(9), which is isotopic to
$K(-2/7,\, 1/2,\,-1/6)$ in Figure 3.1(10), and $K(5)$ is the double
branched cover of $L[1]$ in Figure 1(11), isotopic to $K(3/5,\,
-1/3,\, -1/5)$ in Figure 3.1(12).  This completes the proof for the
three Seifert fibered surgeries on $K(-1/2,\ 1/3,\ 2/5)$.

The proofs for cases (4)--(10) are similar.  Surgeries on
$K(-1/2,\ 1/5,\ 2/5)$ and $K(-1/2,\ 1/5,\ 2/7)$ are given in Figure 3.2,
$K(-1/2,\ 1/3,\ 2/7)$ in Figure 3.3, $K(-1/2,\ 1/3,\ 2/9)$ in Figure 3.4,
$K(-1/2,\ 1/3,\ 2/11)$ in Figure 3.5, and $K(-1/3,\ 1/3,\ 1/4)$ and
$K(-1/3,\ 1/3,\ 1/6)$ in Figure 3.6.

For case (13), the knot can be written as $K = K(1/3, -3/5, 1/3)$ and
can be drawn as in Figure 3.7(1).  This is also a strongly
invertible knot, only that the axis is not the one passing through all
3 tangles as in examples above.  The quotients $\bar X$ and $\bar N$
are shown in Figure 3.7(2), which is isotopic to that in Figure
3.7(3).  Using a symmetric spanning surface as above, one can check
that the longitude projects to a curve of slope $-6$ on $\bdd \bar N$,
hence $K(-5)$ is the double branched cover of $L[-1]$ in Figure
3.7(4), which is isotopic to the Montesinos knot $K(2/5, 2/5, -3/4)$
in Figure 3.7(5).  The proof for the case (11) is similar and is shown
in Figures 3.7(6)--(9).  
 
The $(-3,3,5)$ pretzel knot in case (12) does not seem to be strongly
invertible and hence cannot be proved using the method above.
Fortunately this has been done by Mattman, Miyazaki and Motegi.
Figure 3.7 in [MMM] shows that $1$ surgery on the $(-3,3,5)$ pretzel
knot yields the manifold $M(-1/3, -1/4, 3/5)$.  \qed

\bigskip
\leavevmode

\centerline{\epsfbox{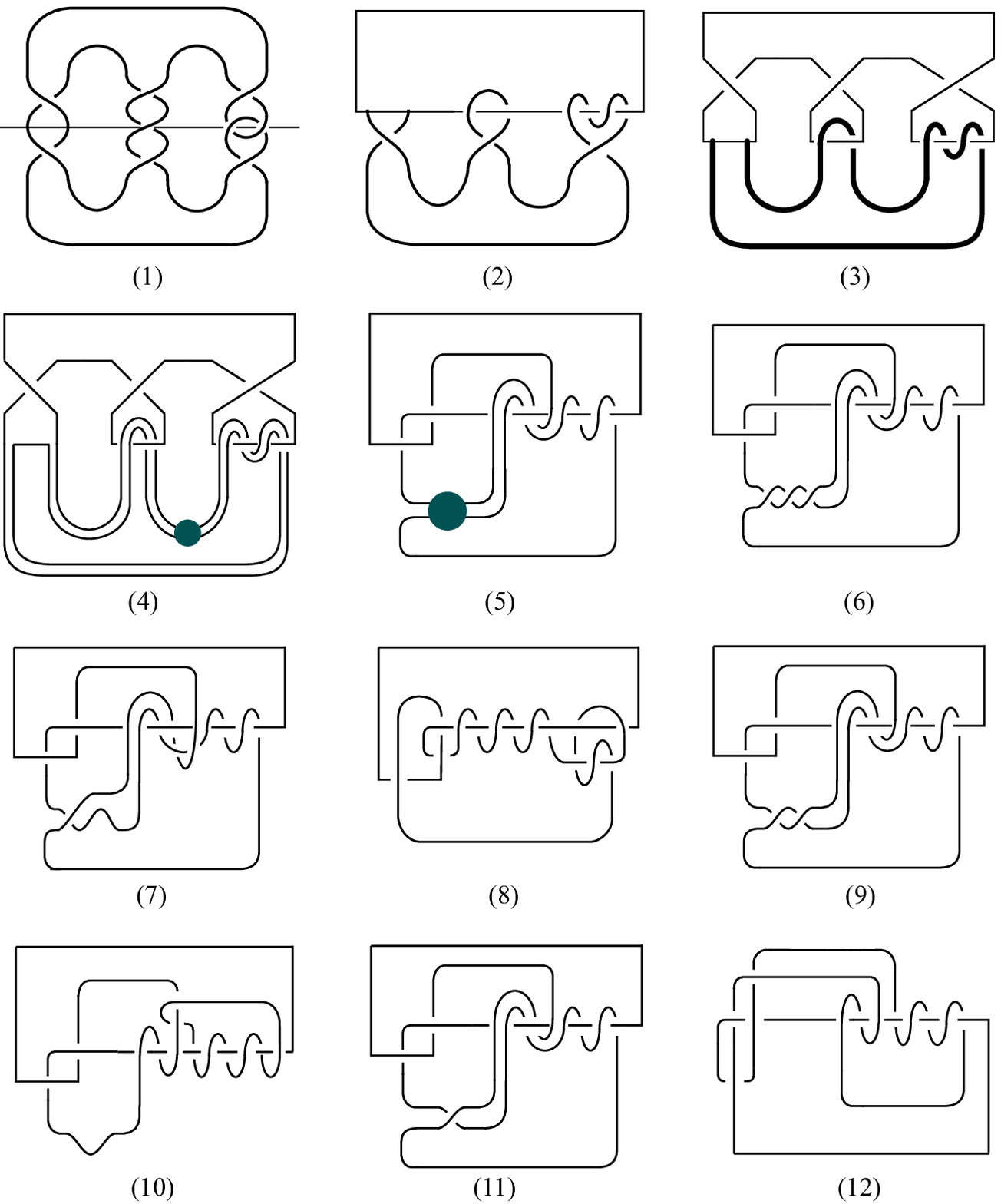}}
\bigskip
\centerline{Figure 3.1}
\bigskip

\bigskip
\leavevmode

\centerline{\epsfbox{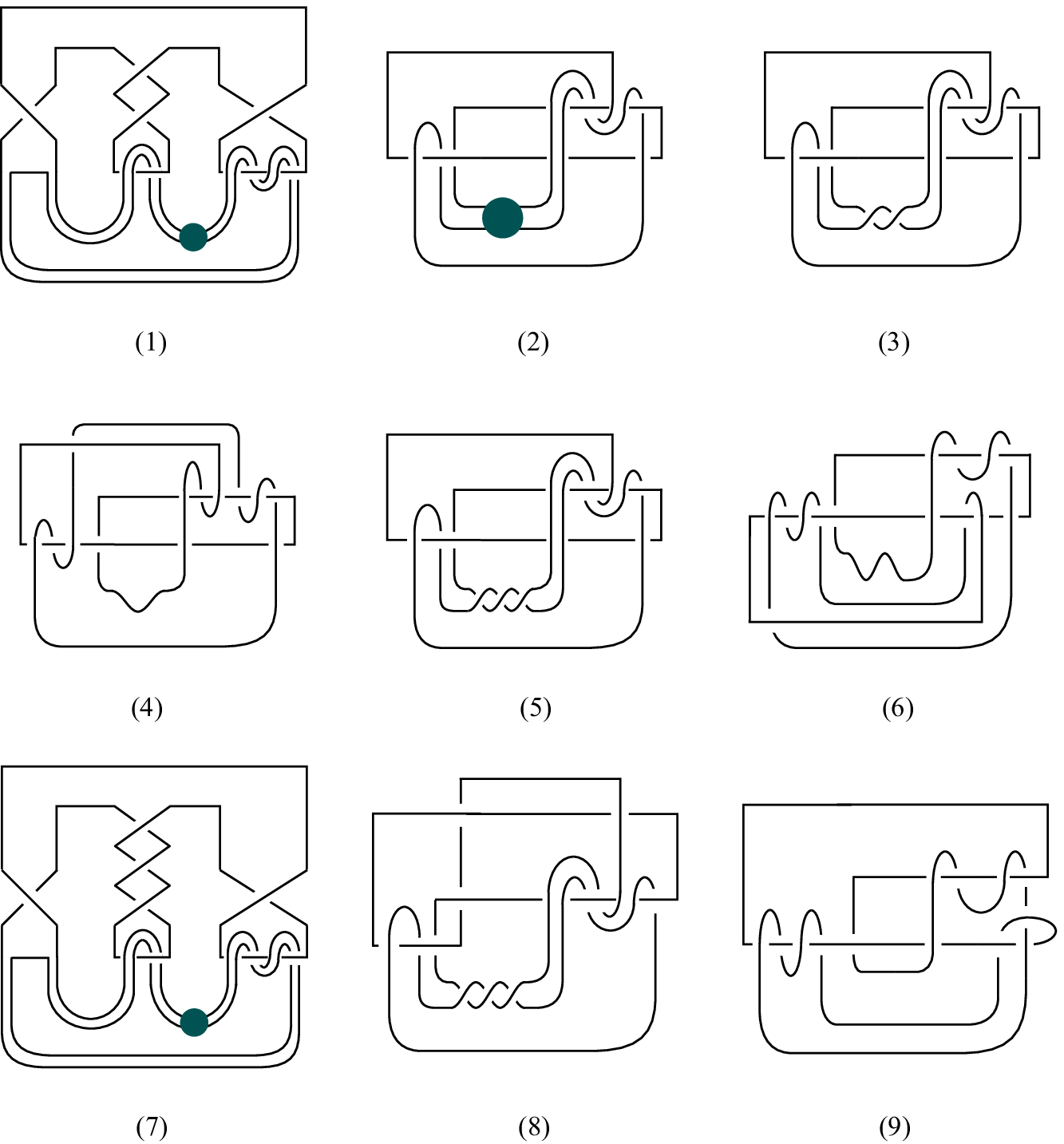}}
\bigskip
\centerline{Figure 3.2}
\bigskip
 
\bigskip
\leavevmode

\centerline{\epsfbox{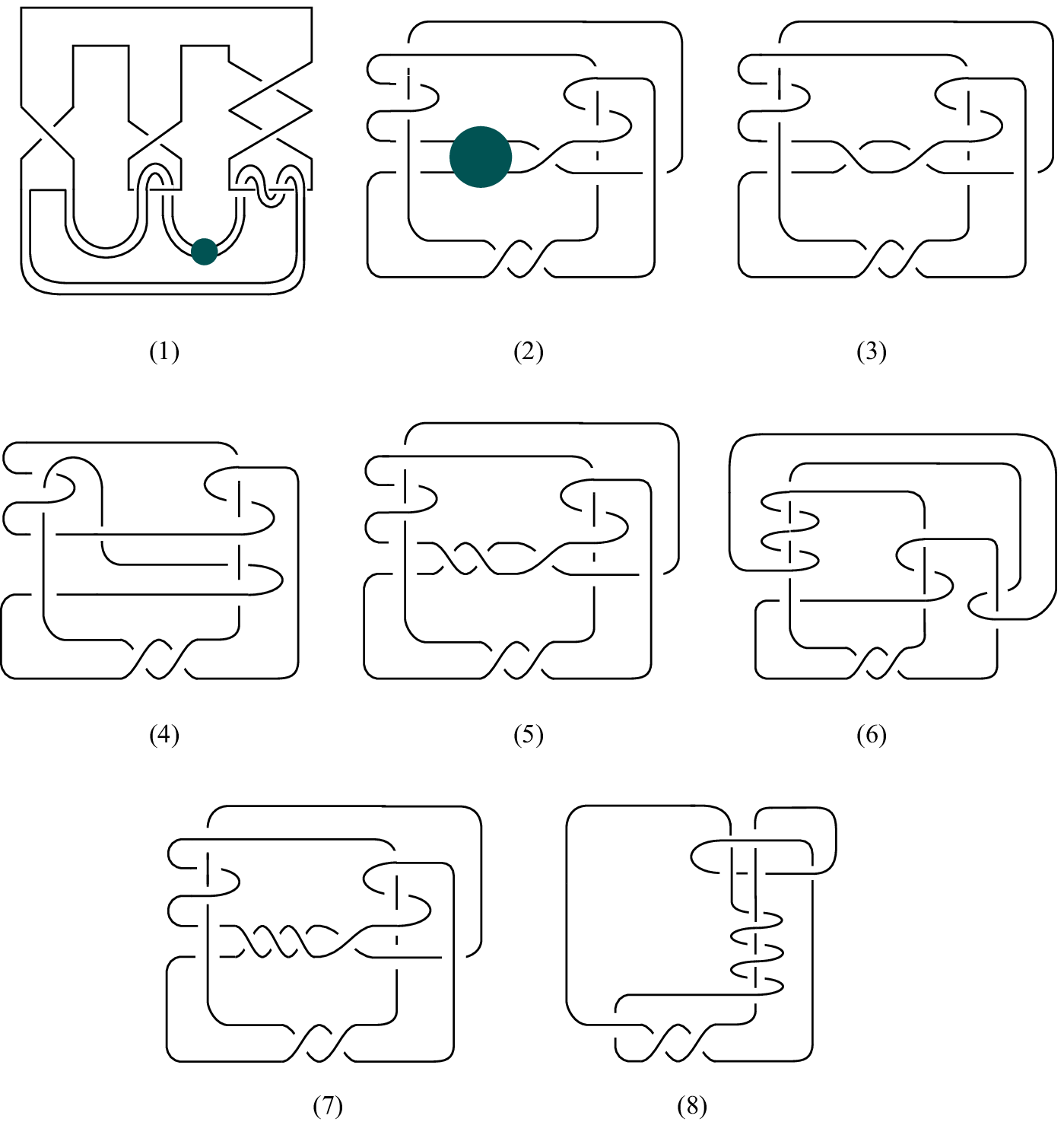}}
\bigskip
\centerline{Figure 3.3}
\bigskip

\bigskip
\leavevmode

\centerline{\epsfbox{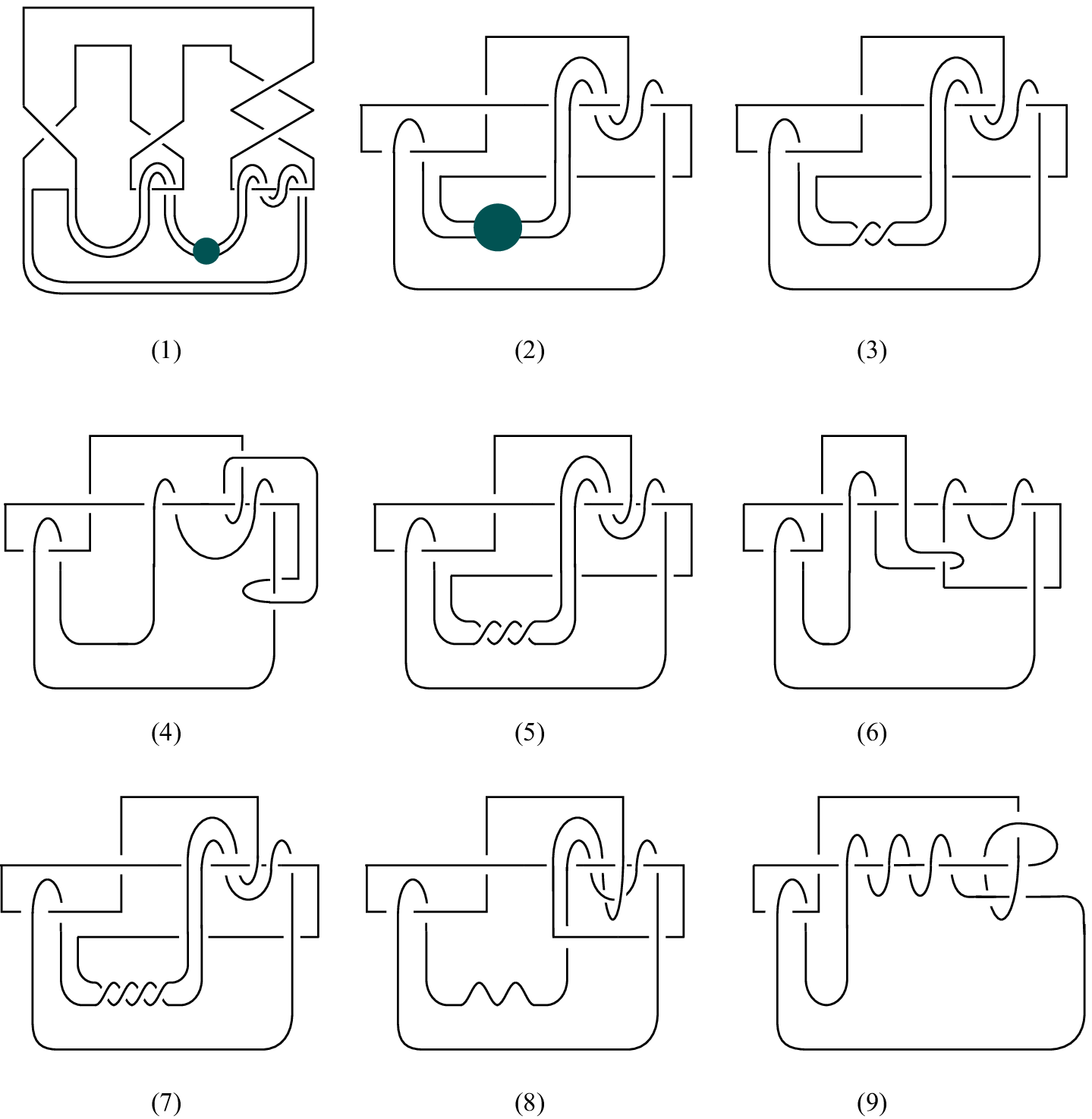}}
\bigskip
\centerline{Figure 3.4}
\bigskip

\bigskip
\leavevmode

\centerline{\epsfbox{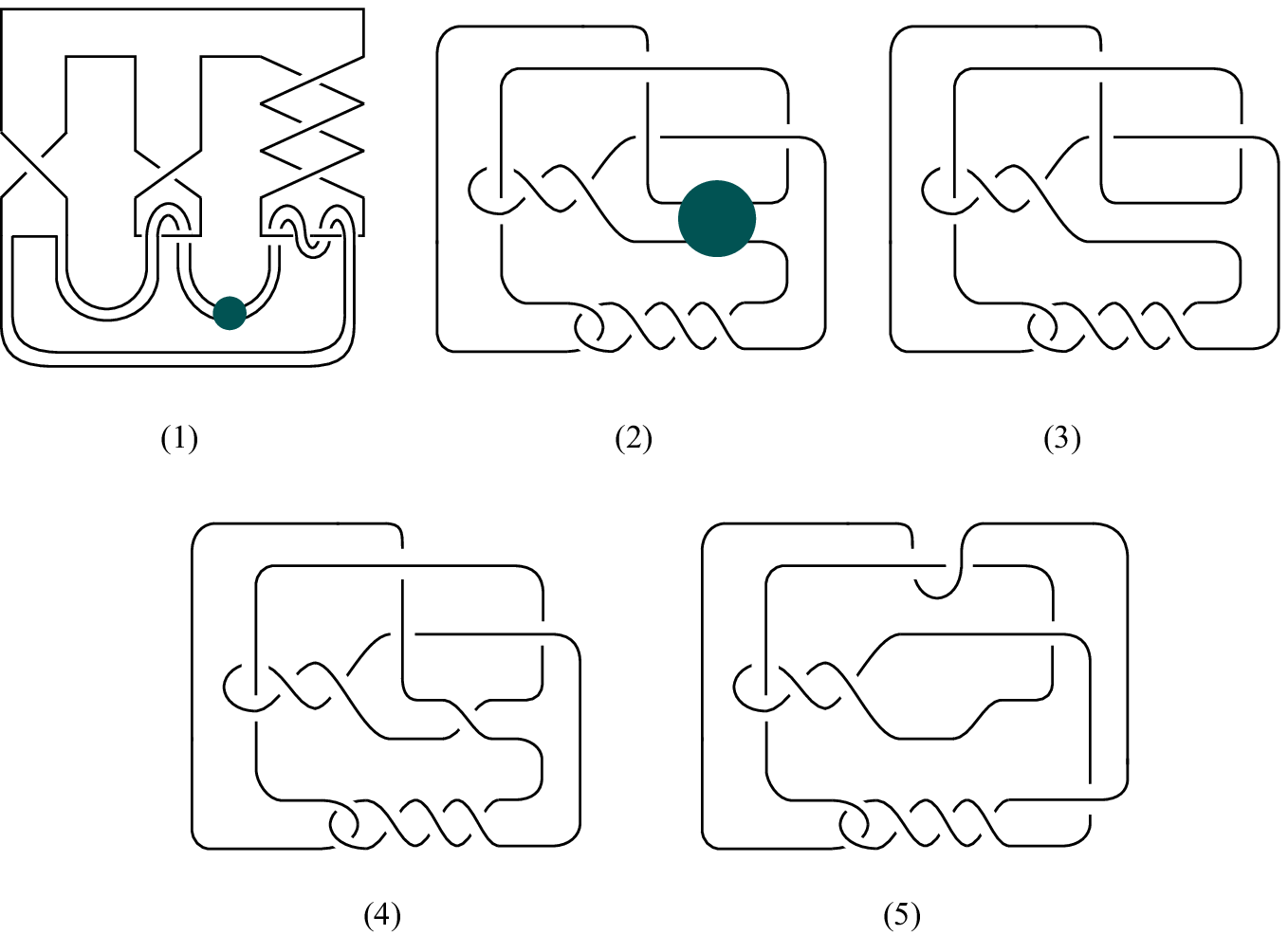}}
\bigskip
\centerline{Figure 3.5}
\bigskip

\bigskip
\leavevmode

\centerline{\epsfbox{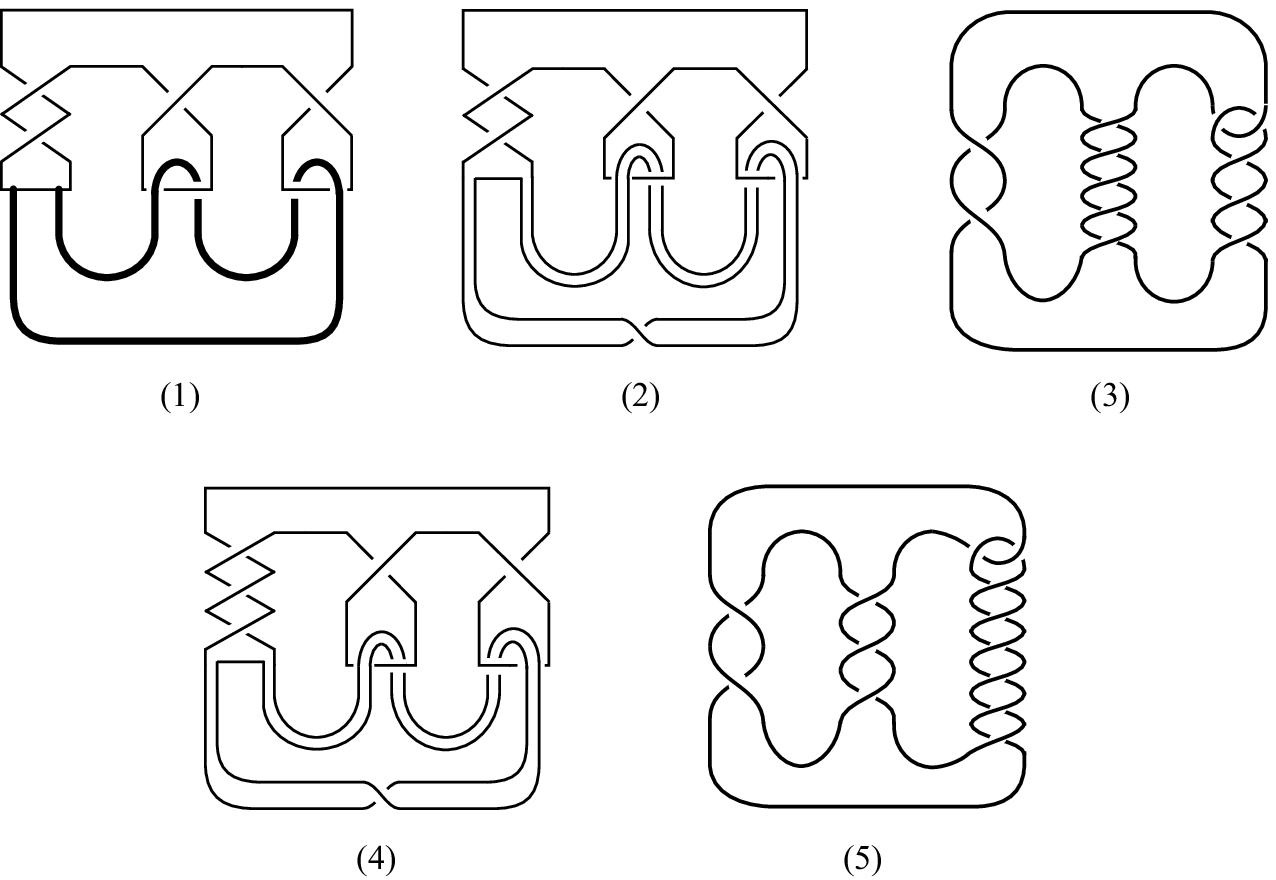}}
\bigskip
\centerline{Figure 3.6}
\bigskip

\bigskip
\leavevmode

\centerline{\epsfbox{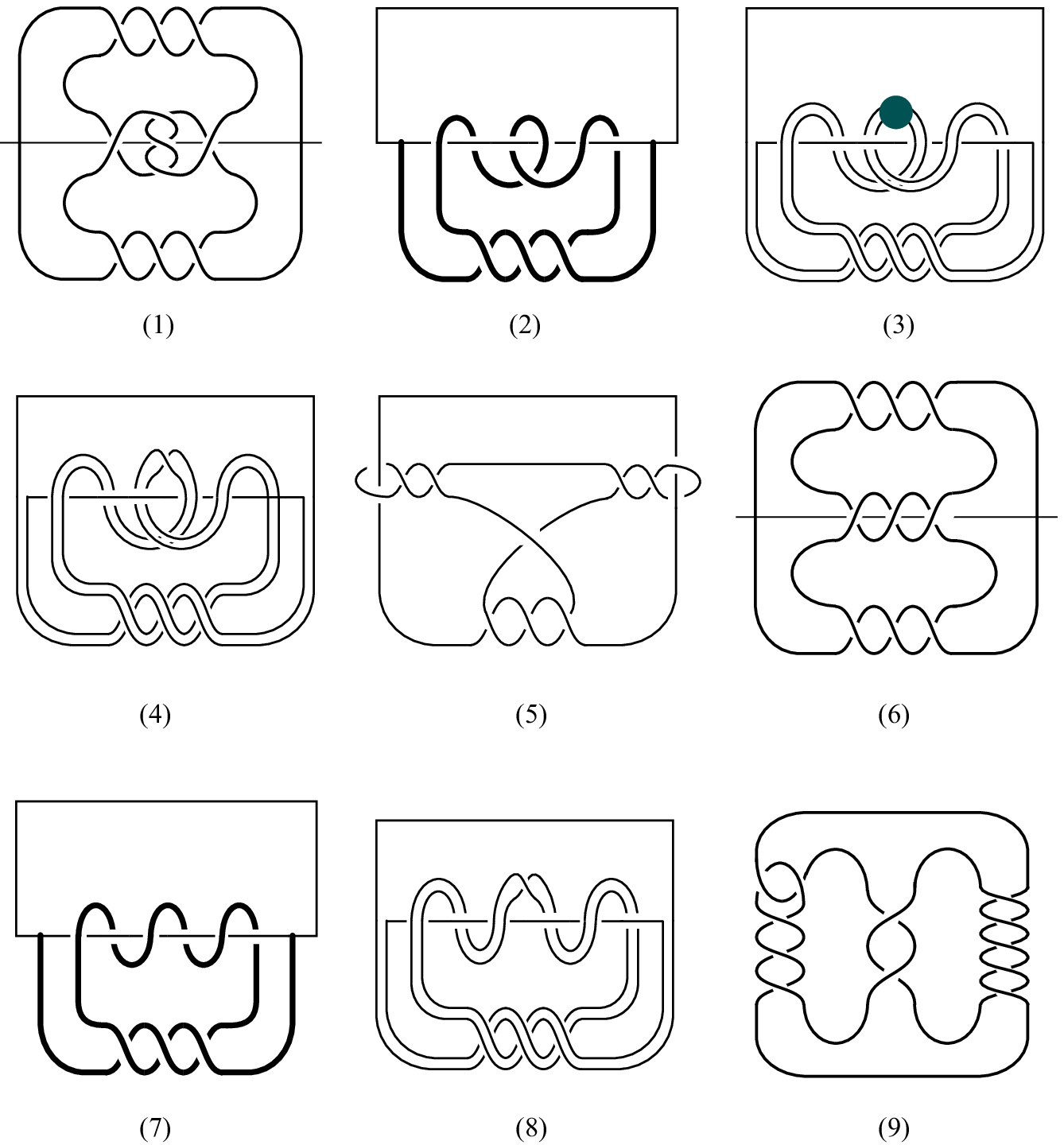}}
\bigskip
\centerline{Figure 3.7}
\bigskip

\section{A conjecture and some computer assistant approach}

Using Snappea or Snappy, one can test the knots in Theorem 2.3 to find
the set of slopes along which Dehn surgeries produce manifolds with
near zero volume.  Snappy volume represents the Gromov norm of the
manifold, hence if the surgery is Seifert fibered then the volume
should be zero.  We have following conjecture.

\begin{conj} A nontrivial Dehn surgery on a hyperbolic Montesinos knot
  of length 3 is Seifert fibered if and only if it is equivalent to
  one of those in table 3.1.
\end{conj}

A computer program {\it Snappex\/} has been written, which combines
{\it Snap\/} of Oliver Goodman [Gm] with a template written by Harriet
Moser [Mos2].  See [Wu7].  {\it Snap\/} uses the {\it Snappea\/} core
of Jeff Weeks [We] and the high precision package {\it Para\/} to
calculate hyperbolic structure for 3-manifolds, while the Moser Script
uses the {\it Snap\/} output as its input and then attempt to verify
the hyperbolicity of the manifold rigorously.  This is based on
[Mos1], in which Moser showed that there is a genuine hyperbolic
structure for the manifold in a neighborhood of the {\it Snap\/}
solution if the latter satisfies certain conditions.  Given a knot
$K$, one can use {\it Snap\/} to find a hyperbolic structure, use
Moser Script to verify it, then use {\it Snap\/} output to find all
slopes of length at most $2\pi$, and then use Moser Script to check
whether each of these is a hyperbolic surgery.  {\it Snappex\/} makes
this procedure automatic.  Thus given a knot $K$, {\it Snappex\/} will
give a list of slopes which contains all possible exceptional slopes.
Assuming the accuracy of the programs involved and the correctness of
compilers, {\it Snappex\/} rigorously proves that Dehn surgery on $K$
along any slope not in the above list must be hyperbolic.

A similar procedure is carried out by {\it Snappex\/} for links of 2
components.  Using this and some theoretical arguments one can show
that $n \leq 9$ for the knots in Theorem 2.3(3).  Consider the link $L
= K' \cup K''$ in Figure 2.1 with $p_1/q_1 = -1/2$ and $p_2/q_2 =
2/5$.  Then the knot in Theorem 2.3(3), which we denote by $K_n$ for
any given $n$, is obtained from $L$ by $-1/n$ surgery on $K''$ Using
Kirby Calculus [Ro] it can be shown that $s$ surgery on $K$ is
equivalent to $(s-4n, -1/n)$ surgery on $L$.  Running {\it Snappex\/}
on this link gives the candidate list $C$, and one can check to see
that if $(r_1, r_2) \in C$ then either $r_2 = 1/n$ with $n \leq 4$, or
$r_1 \in \{-2, -1, 0, 1, 1/0\}$.  We need to show that if $L(r, -1/n)$
is small Seifert fibered for $r=-2,-1,0,1$ then $n \leq 9$.

Consider the case $r=1$.  By Lemma 2.2, the manifold $L(r, \emptyset)$
is hyperbolic for all $r\neq 1/0$; in particular, $M = L(1,
\emptyset)$ is hyperbolic.  Note that $-1$ surgery on the second
component of $L$ yields the knot $K_1 = K(-1/2, 2/5, 1/3)$, which by
Theorem 2.3 has Seifert fibered surgeries of slopes $s = 3,4,5$.  By
the above, we have that $L(-1, -1/1) = K_1(3)$ is non-hyperbolic,
hence by the 8-Theorem $L(-1, -1/n) = K_n(3)$ is hyperbolic for $n>9$.
Similarly for $r = 0, -1$.

The above does not work for $r = -2$.  Fortunately $s = r+4n$ is the
boundary of a non-orientable checkboard spanning surface $F$ with
$\chi(F) = -2$.  Cutting along $F$ produces a handlebody $M$ of genus
3.  Considering $\bdd N(F) \cap M$ as horizontal surface and $\bdd M
\cap \bdd N(K)$ as vertical surface, we obtain a cusped manifold.  It
can be shown that the horizontal surface is incompressible (at least
for $n>2$) and extends to an incompressible surface in the surgered
manifold $K_n(s)$, and $M$ is not an $I$-bundle.  It follows from [Br]
that $K_n(s)$ cannot be small Seifert fibered.  

\medskip

Back to Conjecture 4.1.  We now have a list of a few hundred knots to
check.  {\it Snappex\/} has a command to find a candidate list of
exceptional slopes for all these knots.  There are several hundred
surgeries that remains on this list, which need to be verified using
some other methods.  A few non-integral slopes can be excluded using
Lemma 3.1.  All but a couple of the remaining slopes are integral
slopes, which are shown to be ``apparently hyperbolic'' by Casson's
{\it Geo\/} program [Ca], which provides strong supporting evidence for
the conjecture.

\bigskip

\noindent
Department of Mathematics,  University of Iowa,  Iowa City, IA 52242
\\
Email: {\it wu@math.uiowa.edu}


\begin{thebibliography}{FIKMS} 

{\small

\bibitem[BGZ]{BGZ} S.~Boyer, C.~Gordon and X.~Zhang, {\em Dehn filling
    of knot manifolds containing essential once-punctured tori},
  preprint, ArXiv:1109.5151.

\bibitem[Br]{Br} M.\ Brittenham, {\em Essential laminations in
    Seifert-fibered spaces}, Topology {\bf 32} (1993), 61--85.

\bibitem[BW]{BW} M.\ Brittenham and Y-Q.\ Wu, {\em The classification
    of exceptional Dehn surgeries on 2-bridge knots}, Comm.\ Anal.\
  Geom.\ {\bf 9} (2001), 97--113.

\bibitem [Ca]{Ca} A.\ Casson {\em Geo},
  http://www.math.uiuc.edu/$\sim$nmd/computop/.

\bibitem[CDW]{CDW} M.\ Culler, N.\ Dunfield and J.\ Weeks, {\em
    SnapPy, a computer program for studying the geometry and topology
    of 3-manifolds}, http://snappy.computop.org.

\bibitem[CGLS]{GLS} M.\ Culler, C.\ Gordon, J.\ Luecke and P.\ Shalen,
  {\em Dehn surgery on knots}, Annals Math.\ {\bf 125} (1987),
  237--300

\bibitem[FS]{FS} R.\ Fintushel and R.\ Stern, {\em Constructing lens
    spaces by surgery on knots}, Math.\ Z.\ {\bf 175} (1980), 33--51.

\bibitem[FIKMS]{FIKMS} D.\ Futer, M.\ Ishikawa, Y.\ Kabaya, T.\
  Mattman and K.\ Shimokawa, {\it Finite surgeries on three-tangle
  pretzel knots}, Algebr.\ Geom.\ Topol.\ {\bf 9} (2009), 743--771.

\bibitem[Gm]{Gm} O.\ Goodman, {\em Snap}, http://www.ms.unimelb.edu.au/$\sim$snap/.

\bibitem[HT]{HT} A.\ Hatcher and W.\ Thurston, {\em Incompressible
    surfaces in 2-bridge knot complements}, Invent.\ Math.\ {\bf 79}
  (1985), 225--246.

\bibitem[IJ1]{IJ1} K.\ Ichihara and I.\ Jong, {\it Toroidal Seifert
    fibered surgeries on Montesinos knots}, Comm.\ Anal.\ Geom.\ {\bf
    18} (2010), 579--600.

\bibitem[IJ2]{IJ2} K.\ Ichihara, I.\ Jong, {\em Cyclic and finite
    surgeries on Montesinos knots}, Alg.\ Geom.\ Topol.\ {\bf 9}
  (2009), 731--742.

\bibitem[LM]{LM} M.\ Lackenby and R.\ Meyerhoff, {\em  The maximal
    number of exceptional Dehn surgeries}, Preprint.

\bibitem[MMM]{MMM} T.\ Mattman, K.\ Miyazaki and K.\ Motegi, {\em
    Seifert-fibered surgeries which do not arise from
    primitive/Seifert-fibered constructions}, Trans.\ Amer.\ Math.\
  Soc.\  {\bf 358} (2005), 4045-4055.
 
\bibitem[Mon]{Mon} J.\ Montesinos, {\em Surgery on links and double branched
coverings of $S^3$}, Ann.\ Math.\ Studies {\bf 84} (1975), 227--260.

\bibitem[Mos1]{Mos1} H.\ Moser, {\em Proving a manifold to be hyperbolic once it has
been approximated to be so},  Alg.\ Geom.\ Topol.\ {\bf 9} (2009),
103-133. 

\bibitem[Mos2]{Mos2} ------, {\em Moser Script},
  http://www.math.columbia.edu/$\sim$moser/template.txt.

\bibitem[Ro]{Ro} D. Rolfsen, {\em Knots and Links}, Publish or Perish,
  1990.

\bibitem[We]{We} J.\ Weeks, {\em SnapPea},
  http://www.geometrygames.org/SnapPea/index.html.

\bibitem[Wu1]{Wu1} Y-Q.\ Wu, {\em Dehn surgery on arborescent knots},
  J.\ Diff.\ Geom.\ {\bf 42} (1996), 171--197.

\bibitem[Wu2]{Wu2} ------, {\em Exceptional Dehn surgery on large
    arborescent knots}, Pac.\ J.\ Math.\ {\bf 252} (2011), 219--243.  

\bibitem[Wu3]{Wu3} ------, {\em The classification of toroidal Dehn
    surgeries on Montesinos knots}, Comm.\ Anal.\ Geom.\ {\bf 19}
  (2011), 305-345.

\bibitem[Wu4]{Wu4} ------, {\em Immersed surfaces and Seifert fibered
    surgery on Montesinos knots}, Trans.\ Amer.\ Math.\ Soc.\ (to appear). 

\bibitem[wu5]{Wu5} ------, {\em Persistently laminar branched
    surfaces}, Comm.\ Anal.\ Geom.\ (to appear).

\bibitem[Wu6]{Wu6} ------, {\em Dehn surgery on knots of wrapping
    number 2}, preprint.

\bibitem[Wu7]{Wu7} ------, {\em Snappex}, http://www.math.uiowa.edu/$\sim$wu/snappex/snappex.html.

}

\end{thebibliography}
\enddocument